\documentclass[10pt]{article}

\usepackage{amssymb}
\usepackage{amsmath}
\numberwithin{equation}{section}
\newcommand{\DDD}{\mathcal{D}}

\newcommand{\la}{\lambda}
\newcommand{\va}{\varphi}
\newcommand{\ppp}{\partial}

\newcommand{\hhalf}{\frac{1}{2}}

\newcommand{\weightto}{e^{2s\alpha(x,t_0)}}
\newcommand{\ceta}{\Vert \eta\Vert_{C(\overline{\Omega})}}
\newcommand{\ddd}{\mbox{div}\thinspace}

\newcommand{\weight}{e^{2s\alpha}}

\newcommand{\sumij}{\sum_{i,j=1}^d}

\newcommand{\R}{\mathbb{R}}

\newcommand{\www}{\widetilde}

\newcommand{\ooo}{\overline}
\newcommand{\OOO}{\Omega}

\allowdisplaybreaks
\title
{\bf 
Global Lipschitz stability for an inverse coefficient problem for 
a mean field game system}

\author{$^1$ Oleg Imanuvilov and  
$^2$ Masahiro Yamamoto}

\date{}
\begin{document}
\maketitle
\thanks{
$^1$ Department of Mathematics, Colorado State University, 101 Weber Building, 
Fort Collins CO 80523-1874, USA 
e-mail: {\tt oleg@math.colostate.edu}
\\
%$^2$ Department of Mathematics, City University of Hong Kong, Kowloon, 
%Hong Kong SAR, China email: {\tt hongyliu@cityu.edu.hk}
%\\
$^2$ Graduate School of Mathematical Sciences, The University
of Tokyo, Komaba, Meguro, Tokyo 153-8914, Japan 
e-mail: {\tt myama@ms.u-tokyo.ac.jp}
}
\begin{abstract} 
For an inverse coefficient problem of determining 
a state-varying factor in the corresponding Hamiltonian 
for a mean field game system, 
we prove the global Lipschitz stability by 
spatial data of one component and interior data in an arbitrarily chosen 
subdomain over a time interval.
The proof is based on Carleman estimates with different norms.
\end{abstract} 
\baselineskip 18pt

\section{Introduction}

Let $\OOO \subset \R^d$ be a bounded domain with smooth boundary 
$\ppp\OOO$, and let $T>0$, and $Q:= \OOO\times (0,T)$.
In  article  concerned   a system of
the mean field game:
$$
\left\{ \begin{array}{rl}
& \ppp_tu(x,t) + \Delta u(x,t) - \hhalf p(x)\vert \nabla u(x,t)\vert^2
- c_0(x)v(x,t) = 0 \quad 
\mbox{in $Q$}, \\
& \ppp_tv(x,t) - \Delta v(x,t) 
- \ddd (p(x)v(x,t)\nabla u(x,t)) = 0 \quad 
\mbox{in $Q$},\\ \quad 
&u\vert_{\partial\Omega\times (0,T)} = v
\vert_{\partial\Omega\times (0,T)} = 0.  
\end{array}\right.
                                       \eqno{(1.1)}
$$
In (1.1), $x$ and $t$ are the state and the time variables, 
and $u$ and $v$ denote the value of the game and the population 
density of players respectively (e.g., \cite{ACDPS}, \cite{LL}).
We note that $p(x)$ specifies the Hamiltonian for (1.1).

We consider
\\
{\bf Inverse coefficient problem.}
{\it Let $\omega \subset \OOO$ be an arbitrarily chosen subdomain and 
$t_0 \in (0,T)$, $\delta_0 > 0$ be small.  
One need to determine $p(x)$, $x\in \OOO$ by 
$u\vert_{\omega \times (t_0-\delta_0\, t_0+\delta_0)}$, 
$v\vert_{\omega \times (t_0-\delta_0\, t_0+\delta_0)}$ and
$u(\cdot,t_0)$ in $\OOO$.
}

Henceforth we set $\ppp_i:= \frac{\ppp}{\ppp x_i}$ for $1\le i\le d$,
and define 
$H^{2,1}(Q):= \{ u\in L^2(Q);\, \nabla u, \ppp_i\ppp_ju, \ppp_tu  
\in L^2(Q), \, \, 1\le i,j \le d \}$, and $W^{1,\infty}(\OOO)$
 denotes the Sobolev space of functions whose first partial derivatives are in  $L^\infty(\Omega).$

We state our main result.
\\
{\bf Theorem 1.}
{\it 
For $\ell=1,2$, let $(u_{\ell}, v_{\ell}) \in (H^{2,1}(Q))^2$ satisfy (1.1) 
with $p:= p_{\ell}$.
We assume that 
$$
\left\{ \begin{array}{rl}
& \Vert p_{\ell}\Vert_{W^{1,\infty}(\OOO)} \le M, \quad
\Vert c_0\Vert_{L^{\infty}(\OOO)} \le M, \\
  &\Vert u_{\ell}\Vert_{W^{1,\infty}(Q)}+\Vert\partial_tu_{\ell}\Vert_{L^\infty(0,T;W^{1,\infty}(\Omega))}\le M, \,
\Vert v_{\ell}\Vert_{W^{1,\infty}(Q)}+\Vert\partial_tv_{\ell}\Vert_{L^\infty(0,T;W^{1,\infty}(\Omega))} \le M, \\
& (\ppp_t^ku_{\ell},\ppp_t^kv_{\ell}) \in (H^{2,1}(Q))^2, \quad k=0,1,\,
\ell=1,2
\end{array}\right.
                                \eqno{(1.2)}
$$
and there exists a constant $\delta>0$ such that 
$$
\vert \nabla u_1(x,t_0)\vert \ge \delta, \,\, x\in \OOO 
\quad \mbox{or}\quad 
\vert \nabla u_2(x,t_0)\vert \ge \delta, \,\, x\in \OOO.
$$
Then there exists a constant $C=C(M,\delta) > 0$ such that 
$$
\Vert p_1 - p_2\Vert_{L^2(\OOO)}
\le C(\Vert u_1(\cdot,t_0) - u_2(\cdot,t_0)\Vert_{H^2(\OOO)}
+ \Vert u_1 - u_2\Vert_{H^1(0,T;L^2(\omega))}
+ \Vert v_1 - v_2\Vert_{H^1(0,T;L^2(\omega))}).
$$
}

We emphasize the following features of Theorem 1:
\begin{enumerate}
\item
The observation subdomain $\omega \subset \OOO$ can be 
arbitrarily small.
\item
Lipschitz stability over $\OOO$.
\item
We need the positiveness  of $\vert \nabla u_1\vert$ or 
$\vert \nabla u_2\vert$ only at one moment $t=t_0$.
\item
We do not need spatial data neither $v_1$ nor $v_2$ in $\OOO \times \{t_0\}$.
\end{enumerate}
The features (1) - (3) are inevitable consequences of our methodology (e.g., 
\cite{IY98}).  The last (4) is a new aspect by that 
the unknown is only one spatial function, so that we do not need 
$v(\cdot,t_0)$ but one spatial data are sufficient.   We remark that 
the stability and the uniqueness are open in the cases of $t_0=0$ and $t_0=T$.

As for inverse problems for mean field games, see \cite{ILY1}, \cite{ILY2}, 
\cite{Kl23} - \cite{KLL2},  \cite{LMZ} - \cite{LZ2}.  In particular, 
\cite{Kl2023} proves H\"older stability with extra data compare to the above inverse problem  for the case where
the equation contain special non-local term.

The proof of Theorem 1 is based on Carleman estimates: Lemma 2 
in Section 2 and the modified argument in \cite{IY98}. 

\section{Main Carleman estimate}

Setting $y := u_1 - u_2$, $z:= v_1-v_2$ and $f:= \frac{1}{2}(p_1-p_2)$, 
$g:= \vert \nabla u_1\vert^2$, $h := 2v_1\nabla u_1$,
$r_1:= \frac{1}{2}p_2(\nabla u_1 + \nabla u_2)$,
$r_2:= p_2\nabla u_1$ and $r_3:= p_2v_2.$  From (1.2) we obtain
the  linearized system:
$$
\left\{ \begin{array}{rl}
& \ppp_ty(x,t) + \Delta y(x,t) = c_0z + r_1\cdot \nabla y + g(x,t)f(x) \quad  \mbox{in $Q$},  \\
& \ppp_tz(x,t) - \Delta z(x,t) = \ddd (r_2z+r_3\nabla y) + \ddd (h(x,t)f(x)) 
\quad  \mbox{in $Q$}
\end{array}\right.
                                       \eqno{(2.1)}
$$
with $y = z = 0$ on $\ppp\OOO\times (0,T)$.
By (1.2) we see that 
$$
\ppp_t^jr_k, \ppp_t^jg, \ppp_t^jh \in L^{\infty}(Q) \quad 
\mbox{with $1\le k \le 3$ and $j=0,1$}.               \eqno{(2.2)}
$$
For the Carleman estimates, we introduce $\eta \in C^2(\ooo{\OOO})$ such that 
$\eta > 0$ in $\OOO$, $\eta\vert_{\ppp\OOO} = 0$ and 
$\vert \nabla \eta \vert > 0$ on $\ooo{\OOO \setminus \omega_0}$,
where $\omega_0$ is some subdomain such that 
$\ooo{\omega_0} \subset \omega$.   (For existence of such a function see e.g. \cite{Im}.)
Without loss of generality, we can assume that $t_0 = \hhalf T$ and 
$(t_0-\delta_0\, t_0+\delta_0) = (0,T)$.

Fixing a constant $\la > 0$ sufficiently large, we
set $\mu(t):= t(T-t)$, $Q_{\omega}:= \omega \times (0,T)$, and
$$
\va(x,t):= \frac{e^{\la\eta(x)}}{\mu(t)}, \quad
\alpha(x,t):= \frac{e^{\la\eta(x)} - e^{2\la \ceta}}{\mu(t)}, \quad
(x,t) \in Q.
$$
Then we state Carleman estimates for single parabolic equations, which 
can be both backward and forward.
\\
{\bf Lemma 1.}
{\it 
Let $\widetilde y, \widetilde z \in  H^{2,1}(Q)$ and  
$\widetilde y\vert_{\partial\Omega\times (0,T)} =
\widetilde z\vert_{\partial\Omega\times (0,T)} = 0$.  
There exist constants $s_0>0$ and $C>0$ such that for all $s>s_0$, we have
$$
\int_Q \left(\frac{1}{s\va}\left( \sumij \vert \ppp_i\ppp_j\widetilde y\vert^2
+ \vert \ppp_t\widetilde y\vert^2\right) + s\va\vert \nabla \widetilde y\vert^2
+ s^3\va^3\vert \widetilde y\vert^2\right) \weight dxdt
$$
$$
\le C\int_Q \vert \ppp_t\widetilde y + \Delta\widetilde  y\vert^2 \weight dxdt 
+ C\int_{Q_{\omega}} s^3\va^3 \vert\widetilde y\vert^2 \weight dxdt
                            \eqno{(2.3)}
$$
and 
$$
\int_Q \left( \frac{1}{s\va} \vert \nabla\widetilde  z\vert^2
+ s\va\vert \widetilde z\vert^2\right) \weight dxdt
\le C\int_Q \vert G\vert^2 \weight dxdt 
+ C\int_{Q_{\omega}} s\va \vert \widetilde z\vert^2 \weight dxdt
                                       \eqno{(2.4)}
$$
where $\ppp_t\www{z} - \Delta \www{z} = \ddd G$ in $Q$.
}

The Carleman estimate (2.3) can be proved by 
applying $\alpha(x,t) = \alpha(x,T-t)$ and $\va(x,t) = \va(x,T-t)$ for $(x,t)
\in Q$ to the Carleman estimate in \cite{Im}, 
while the proof of (2.4) is 
found in Imanuvilov and Yamamoto \cite{IY3}.

Henceforth $C>0$ denotes generic constants independent of $s>0$.
We define 
$D(y,z):= \int_{Q_{\omega}} (s^3\va^3\vert y\vert^2
+ s\va\vert z\vert^2) \weight dxdt$.
Setting $y_1:= \ppp_ty$ and $z_1:= \ppp_tz$, and  differentiating equations (2.1) respect to $t$ we have
$$
\left\{ \begin{array}{rl}
& \ppp_ty_1 + \Delta y_1 = c_0z_1 
+ r_1\cdot\nabla y_1 + (\ppp_tr_1)\cdot \nabla y 
+ (\ppp_tg)f \quad \mbox{in $Q$}, \\
& \ppp_tz_1 - \Delta z_1 = \ddd (r_2z_1+r_3\nabla y_1) 
+ \ddd ((\ppp_tr_2)z + (\ppp_tr_3)\nabla y) + \ddd((\ppp_th)f)
\quad \mbox{in $Q$}.
\end{array}\right.
                            \eqno{(2.5)}
$$
Applying (2.3) and (2.4) to the first and the second 
equations in (2.1) respectively, in terms of (2.2) we have
$$
\int_Q \left(\frac{1}{s\va} \vert \ppp_ty\vert^2
 + s\va\vert \nabla y\vert^2
+ s^3\va^3\vert y\vert^2\right) \weight dxdt
\le C\int_Q \vert z\vert^2 \weight dxdt
+ C\int_Q \vert f\vert^2 \weight dxdt + CD(y,z)  \eqno{(2.6)}
$$
and
$$
\int_Q s\va\vert z\vert^2 \weight dxdt
\le C \int_Q (\vert y\vert^2 + \vert \nabla y\vert^2) \weight dxdt 
+ C\int_Q \vert f\vert^2 \weight dxdt + CD(y,z).     \eqno{(2.7)}
$$
Here in terms of (2.2), we estimate $\int_Q \vert r_1\cdot \nabla y\vert^2
\weight dxdt \le C\int_Q \vert \nabla y\vert^2\weight dxdt$, and this
term can be absorbed into $\int_Q s\va\vert \nabla y\vert^2 \weight dxdt$ 
on the left-hand side of (2.3) for large $s>0$. 
Throughout the proof, we repeat similar 
estimation with absorption thanks to the large parameter $s>0$.

Adding (2.6) and (2.7), and choosing $s>0$ large, we can absorb the 
resulting term on the right-hand side
into the left-hand side, so that 
$$
\int_Q \left(\frac{1}{s\va} \vert \ppp_ty\vert^2
 + s\va\vert \nabla y\vert^2
+ s^3\va^3\vert y\vert^2 + s\va \vert z\vert^2\right) \weight dxdt
\le C\int_Q \vert f\vert^2 \weight dxdt + CD(y,z).    \eqno{(2.8)}
$$
Next, the application of Lemma 1 to (2.5) yields
$$
 \int_Q \left( \frac{1}{s\va} \vert \ppp_ty_1\vert^2
 + s\va\vert \nabla y_1\vert^2
+ s^3\va^3\vert y_1\vert^2\right) \weight dxdt
$$
$$
\le C\int_Q (\vert \nabla y\vert^2 + \vert z_1\vert^2)
 \weight dxdt
+ C\int_Q \vert f\vert^2 \weight dxdt + CD(y_1,z_1)  \eqno{(2.9)}
$$
and
$$
\int_Q s\va\vert z_1\vert^2 \weight dxdt
\le C \int_Q (\vert z\vert^2 + \vert \nabla y\vert^2 + \vert \nabla y_1\vert^2) \weight dxdt 
+ C\int_Q \vert f\vert^2 \weight dxdt + CD(y_1,z_1).     \eqno{(2.10)}
$$
Adding (2.9) and (2.10), we can absorb the terms $\vert z_1\vert^2$,
$\vert \nabla y_1\vert^2$ on the right-hand side
into the left-hand side, we can obtain
\begin{align*}
& \int_Q \left( \frac{1}{s\va} \vert \ppp_ty_1\vert^2
 + s\va\vert \nabla y_1\vert^2
+ s^3\va^3\vert y_1\vert^2 + s\va\vert z_1\vert^2\right) \weight dxdt\\
\le& C\int_Q (\vert \nabla y\vert^2 + \vert z\vert^2)
 \weight dxdt
+ C\int_Q \vert f\vert^2 \weight dxdt + CD(y_1,z_1).
\end{align*}
Substituting (2.8) into the first term on the right-hand side, we reach
\\
{\bf Lemma 2 (key Carleman estimate.}
{\it 
There exist constants $s_0> 0$ and $C>$ such that 
$$
 \int_Q \biggl( \frac{1}{s\va} \vert \ppp_t^2y\vert^2
 + s^3\va^3\vert \ppp_ty\vert^2)
+ s\va(\vert z\vert^2 + \vert \ppp_tz\vert^2) \biggr)
\weight dxdt
\le C\int_Q \vert f\vert^2 \weight dxdt + C(D(y,z) + D(y_1,z_1))
$$
for all $s>s_0$.
}

\section{Completion of the proof of Theorem 1.}

By $e^{2s\alpha(x,0)} = 0$ for $x\in \OOO$, and 
$\vert \ppp_t\va\vert\le C\va^2, \vert \ppp_t\alpha \vert \le C\va^2$ in $Q$, we have
\begin{align*}
& \int_{\OOO} \va(x,t_0)^{-1} \vert \ppp_ty(x,t_0)\vert^2 e^{2s\alpha(x,t_0)}
dx
= \int^{t_0}_0 \frac{d}{dt} \left( \int_{\OOO} \va^{-1}\vert \ppp_ty\vert^2
\weight dx\right) dt\\
=& \int^{t_0}_0\int_{\OOO} (-(\ppp_t\va)\va^{-2}\vert \ppp_ty\vert^2
+2s \va^{-1}\vert \ppp_ty\vert^2(\ppp_t\alpha)
+ 2\va^{-1}(\ppp_ty)(\ppp_t^2y)) \weight dxdt\\
\le& C\int_Q (\vert \ppp_ty\vert^2 + s\va\vert \ppp_ty\vert^2 
+ \vert \ppp_ty\vert \vert \ppp_t^2y\vert) \weight dxdt
\le C\int_Q (s\va\vert \ppp_ty\vert^2 
+ \vert \ppp_ty\vert \vert \ppp_t^2y\vert) \weight dxdt\\
\le& C\int_Q (s\va\vert \ppp_ty\vert^2 + \frac{1}{s\va} \vert \ppp_t^2y\vert^2)
\weight dxdt.
\end{align*}
Here we used
$\vert \ppp_ty\vert\vert \ppp_t^2y\vert 
= \left(\frac{1}{\sqrt{s\va}}\vert \ppp_t^2y\vert\right)
( \sqrt{s\va}\vert \ppp_ty\vert)
\le \frac{1}{2}\left( \frac{1}{s\va}\vert \ppp_t^2y\vert^2
+ s\va\vert \ppp_ty\vert^2\right)$.
\\
Therefore, by $\min_{x\in \ooo{\OOO}} \va(x,t_0)^{-1}>0$,
Lemma 2 yields
$$
\int_{\OOO} \vert \ppp_ty(x,t_0)\vert^2 e^{2s\alpha(x,t_0)} dx
\le C\int_{\OOO} \va(x,t_0)^{-1} \vert \ppp_ty(x,t_0)\vert^2 
e^{2s\alpha(x,t_0)} dx
\le C\int_Q \vert f\vert^2 \weight dxdt + C_s\mathcal{D}.
$$
Here and henceforth we set 
$\DDD:= \Vert y\Vert^2_{H^1(0,T;L^2(\omega))}
+ \Vert z\Vert^2_{H^1(0,T;L^2(\omega))}$.

We can assume that $\vert \nabla u_1(x,t_0)\vert > 0$ for 
$x\in \ooo{\OOO}$.  Then, since 
$\vert g(x,t_0)\vert = \vert \nabla u_1(x,t_0)\vert^2 > 0$ 
for $x \in \ooo{\OOO}$ and 
$g(x,t_0)f(x) = \ppp_ty(x,t_0) + (\Delta y - r_1\nabla y)(x,t_0)
- c_0(x)z(x,t_0)$ for $x\in \OOO$,
we obtain
$$
\int_{\OOO} \vert f(x) \vert^2 e^{2s\alpha(x,t_0)} dx
\le C\int_Q \vert f\vert^2 \weight dxdt + C_s\mathcal{D}
+ C_s\Vert y(\cdot,t_0)\Vert^2_{H^2(\OOO)}
+ C\int_{\OOO} \vert z(x,t_0)\vert^2 \weightto dx.
                                          \eqno{(3.1)}
$$
Next 
\begin{align*}
& \int_{\OOO} \vert z(x,t_0)\vert^2 e^{2s\alpha(x,t_0)}dx
\le C \int_{\OOO} \va(x,t_0)^{-1} \vert z(x,t_0)\vert^2 e^{2s\alpha(x,t_0)}
dx\\
= & C\int^{t_0}_0 \frac{d}{dt} \left( \int_{\OOO} \va^{-1}\vert z\vert^2
\weight dx\right) dt\\
=& C\int^{t_0}_0\int_{\OOO} (-(\ppp_t\va)\va^{-2}\vert z\vert^2
+ \va^{-1}\vert z\vert^2 2s(\ppp_t\alpha)
+ \va^{-1}2z(\ppp_tz)) \weight dxdt\\
\le& C\int_Q (s\va\vert z\vert^2 + \vert z\vert \vert \ppp_tz\vert)\weight dxdt
\le C\int_Q (s\va\vert z\vert^2 + \vert z\vert^2
+ \vert \ppp_tz\vert^2) \weight dxdt.
\end{align*}
Therefore Lemma 2 yields 
$\int_{\OOO} \vert z(x,t_0)\vert^2 \weightto dx 
\le C\int_Q \vert f\vert^2 \weight dxdt + C_s\DDD$,
with which (3.1) implies
$$
\int_{\OOO} \vert f(x)\vert^2 \weightto dx
\le C\int_Q \vert f\vert^2 \weight dxdt 
+ C_s(\DDD + \Vert y(\cdot,t_0)\Vert^2_{H^2(\OOO)})
                                  \eqno{(3.2)}
$$
for all large $s>0$.

On the other hand,
$$
\int_Q \vert f(x)\vert^2 e^{2s\alpha(x,t)} dxdt 
= \int_{\OOO} \vert f(x)\vert^2 \weightto 
\left( \int^T_0 e^{2s(\alpha(x,t)-\alpha(x,t_0))} dt
\right) dx.
$$
Since $\mu(t_0) > \mu(t)$ for $t\ne t_0$, we verify
$\alpha(x,t_0) - \alpha(x,t)
\ge C_0\left( \frac{1}{\mu(t)} - \frac{1}{\mu(t_0)}\right)$
for $(x,t) \in Q$, where $C_0:= e^{2\la\ceta} - e^{\la\ceta}$.

Hence,
$$
\int^T_0 e^{2s(\alpha(x,t)-\alpha(x,t_0))} dt
\le \int^T_0 \exp\left( -2sC_0\left( \frac{1}{\mu(t)} - \frac{1}{\mu(t_0)}
\right)\right) dt
$$
for $x\in \OOO$.
Since 
$\lim_{s\to \infty} 
\exp\left( -2sC_0\left( \frac{1}{\mu(t)} - \frac{1}{\mu(t_0)}\right)\right)
= 0$ if $t\ne t_0$ and 
$\exp\left( -2sC_0\left( \frac{1}{\mu(t)} - \frac{1}{\mu(t_0)}\right)\right)
\le 1$ for $s>0$ and $0 \le t \le T$, the Lebesgue convergence theorem yields
$\sup_{x\in \OOO} \int^T_0 e^{2s(\alpha(x,t) - \alpha(x,t_0))} dt 
= o(1)$ as $s\to \infty$,
and so 
$\int_Q \vert f\vert^2 \weight dxdt 
= o(1)\int_{\OOO} \vert f\vert^2 \weightto dx$.
Substituting this into (3.2) and choosing $s>0$ large, we reach 
$\int_{\OOO} \vert f\vert^2 \weightto dx
\le C_s(\DDD + \Vert y(\cdot,t_0)\Vert^2_{H^2(\OOO)})$.
$\blacksquare$

{\bf Acknowledgments.}
The work was supported by Grant-in-Aid for Scientific Research (A) 20H00117 
and Grant-in-Aid for Challenging Research (Pioneering) 21K18142 of 
Japan Society for the Promotion of Science.
%The authors thank the anonymous referees for careful reading and 
%valuable comments.

%

\begin{thebibliography}{99} %

\bibitem{ACDPS}
Y. Achdou, P. Cardaliaguet, F. Delarue, A. Porretta and F. Santambrogio, 
{\it Mean field games},
Cetraro, Italy 2019, Lecture Notes in Mathematics, C.I.M.E. 
Foundation Subseries, Volume 2281, Springer, 2019.

\bibitem{Im}
O.Y. Imanuvilov, {\it Controllability of parabolic equations,} Sbornik Math. 
{\bf 186} (1995) 879-900.

\bibitem{ILY1}
O.Y. Imanuvilov, H. Liu and M. Yamamoto, 
{\it Unique continuation for a mean field game system},
to appear in Applied Math. Letters.
https://doi.org/10.1016/j.aml.2023.108757

\bibitem{ILY2}
O.Y. Imanuvilov, H. Liu and M. Yamamoto, 
{\it Lipschitz stability for determination pf states and 
inverse source problem for the the mean field game equations},
preprint, 2023.
 
\bibitem{IY3}
O.Y. Imanuvilov and M. Yamamoto,
{\it Carleman inequalities for parabolic equations in Sobolev spaces of 
negative order and exact controllability for semilinear parabolic equations}, 
Publ. Res. Inst. Math. Sci. {\bf 39} (2003) 227-274.

\bibitem{IY98}
O.Y. Imanuvilov and M. Yamamoto,
{\it Lipschitz stability in inverse parabolic problems by the Carleman
estimate}, Inverse Problems {\bf 14} (1998) 1229-1245.

\bibitem{Kl23}
M. V. Klibanov, {\it 
The mean field games system: Carleman estimates, Lipschitz stability and 
uniqueness,} published online in J. Inverse and Ill-posed Problems, 2023 
arXiv:2303.03928

\bibitem{Kl2023}
K.V. Klibanov, {\it
A coefficient inverse problem for the mean field games system},
preprint, arXiv:2306.03349

\bibitem{KlAv}
M. V. Klibanov and Y. Averboukh,
{\it Lipschitz stability estimate and uniqueness in the retrospective analysis 
for the mean field games system via two Carleman estimates,}
preprint arXiv:2302.10709

\bibitem{KLL1}
M. V. Klibanov, J. Li and H. Liu, {\it 
H\"older stability and uniqueness for the mean field games system 
via Carleman estimates,} preprint, arXiv:2304.00646

\bibitem{KLL2}
M. V. Klibanov, J. Li and H. Liu,
{\it On the mean field games system with the lateral Cauchy data via 
Carleman estimates,} preprint, arXiv:2303.07556 

\bibitem{LL}
J.-M. Lasry and P.-L. Lions,{\it  Mean field games}, 
Japanese Journal of Mathematics, {\bf 2} (2007) 229-260.

\bibitem{LMZ}
H. Liu, C. Mou and S. Zhang, {\it Inverse problems for mean field games, }
preprint, arXiv:2205.11350

\bibitem{LY}
H. Liu and M. Yamamoto, 
{\it Stability in determination of states for the 
mean field game equations,} preprint, arXiv:2304.05896
to appear in Communications on Analysis and Computation
doi:10.3934/cac.2023009

\bibitem{LZ1}
H. Liu and S. Zhang, {\it On an inverse boundary problem for mean field games,}
preprint, arXiv:2212.09110 

\bibitem{LZ2}
H. Liu and S. Zhang, {\it Simultaneously recovering running cost and 
Hamiltonian 
in mean field games system,} preprint,
arXiv:2303.13096

%\bibitem{Y09}
%M. Yamamoto, {\it Carleman estimates for parabolic equations and 
%applications,} Inverse Problems {\bf 25} (2009) 123013.

\end{thebibliography}
\end{document}